\documentclass[a4paper,10pt]{article}
\usepackage[english]{babel}
\usepackage[latin1]{inputenc}
\usepackage{amsfonts}
\usepackage{amsmath}
\usepackage{amsthm}
\usepackage{amssymb}
\usepackage{latexsym}
\usepackage[dvips]{graphicx}
\usepackage{graphics}
\usepackage{psfig}
\usepackage{epsfig}
\usepackage{natbib}
\newtheorem{defn}{Definition}[section]
\newtheorem{lem}{Lemma}[section]
\newtheorem{thm}{Theorem}[section]
\newtheorem{prop}{Proposition}[section]

\newcommand{\et}{{\it et al.\ }}
\newcommand{\Hbf}{{\bf H}}
\newcommand{\wpi}{\widehat{\pi}}
\newcommand{\wL}{\widehat{L}}
\newcommand{\wR}{\widehat{R}}
\newcommand{\wG}{\widehat{G}}
\newcommand{\womega}{\widehat{\omega}}

\newcommand{\E}{\mathbb{E}}

\newcommand{\Abf}{{\bf A}}
\newcommand{\lnorm}{||}
\renewcommand{\wp}{\widehat{p}}
\renewcommand{\P}{\mathbb{P}}
\renewcommand{\wpi}{\widehat{\pi}}
\def\ind{\hbox{ 1\hskip -2.7pt {I}}}
\def\lbr{\hbox{ $\mathrm{[}$\hskip -4.2pt {$\mathrm{[}\,$}}}
\def\rbr{\hbox{$\,\mathrm{]}$\hskip -4.2pt {$\mathrm{]}\,$}}}
\def\argmin{\mathop{\mathrm{Argmin}}}
\begin{document}
\begin{center}
{\Large
    {\sc A leave-p-out based estimation\\ of \\
    the proportion of null hypotheses} }
\bigskip

Alain CELISSE \ \ Stéphane ROBIN\\
\medskip
July 2007

\end{center}
\bigskip
\begin{abstract}
In the multiple testing context, a challenging problem is the
estimation of the proportion $\pi_0$ of true-null hypotheses. A
large number of estimators of this quantity rely on identifiability
assumptions that either appear to be violated on real data, or may
be at least relaxed. Under independence, we propose an estimator
$\widehat{\pi}_0$ based on density estimation using both histograms
and cross-validation.\\
Due to the strong connection between the false discovery rate (FDR)
and $\pi_0$, many multiple testing procedures (MTP) designed to
control the FDR may be improved by introducing an estimator of
$\pi_0$. We provide an example of such an improvement (plug-in MTP)
based on the procedure of Benjamini and Hochberg. Asymptotic
optimality results may be derived for both $\widehat{\pi}_0$ and the
resulting plug-in procedure. The latter ensures the desired
asymptotic control of the FDR, while it is
more powerful than the BH-procedure.\\
Finally, we compare our estimator of $\pi_0$ with other widespread
estimators in a wide range of simulations. We obtain better results
than other tested methods in terms of mean square error (MSE) of the
proposed estimator. Finally, both asymptotic optimality results and
the interest in tightly estimating $\pi_0$ are confirmed
(empirically) by results obtained with the plug-in MTP.
\end{abstract}
%
\section*{Introduction}
Multiple testing problems arise as soon as several hypotheses are
tested simultaneously. Like in test theory, we are concerned with
the control of type-I errors we may commit in falsely rejecting any
tested hypothesis. Post-genomics, astrophysics or neuroimaging are
typical areas in which multiple testing problems are encountered.
For all these domains, the number of tests may be of the order of
several thousands. Suppose we are testing each of $m$ hypotheses at
level $0<\alpha<1,$ the probability of at least one false positive
(\textit{e.g.} false rejection) may equal $m\alpha$ in the worst
case. A possible way to cope with this is to use the Bonferroni
procedure (\cite{DPSB}), which consists in testing each hypothesis
at level $\alpha/m.$ However,
this method is known to be drastically conservative.\\
Since we may be more interested in controlling the proportion of
false positives among rejections rather than the total number of
false positives itself, Benjamini and Hochberg \cite{BH} introduced
the false discovery rate (FDR), defined by
\begin{eqnarray*}
FDR=\mathbb{E}\left[\frac{FP}{1\vee R}\right],
\end{eqnarray*} where $a\vee b=\max(a,b),$ $FP$ denotes the number of false positives and
$R$ is the total number of rejections. A large part of the literature
is devoted to the building of multiple testing procedures (MTP) that
upper bound FDR as tightly as possible (\cite{BKY,BY}). For instance,
that of Benjamini and Hochberg (BH-procedure) \cite{BH} ensures the
following inequality under independence
$$FDR\leq\pi_0 \alpha\leq \alpha,$$ where $\pi_0$ denotes the unknown
proportion of true null hypotheses, while $\alpha$ is the actual
level at which we want to control the FDR. Since $\pi_0$ is unknown,
the BH-procedure suffers some loss in power, which is all the more
deep as $\pi_0$ is small. A natural idea to overcome this drawback
is the computation of an accurate $\pi_0$ estimator, which would be
plugged in the procedure. Thus $\pi_0$ appears as a crucial quantity
that is to be estimated, hence the large amount of existing
estimators. We refer to \cite{LLF,Bro} for reviews on this topic.
The randomness of this estimation needs to be taken into account in
the assessment of the procedure performance (\cite{GW04,Sto02}).\\
In many of quite recent papers about multiple testing (see
\cite{Bro,Efr04,ETST,GW04}), a two-component mixture density is used
to describe the behaviour of p-values associated with the $m$ tested
hypotheses. As usual for mixture models, we need an assumption that
ensures the identifiability of the model parameters. Thus, most of
$\pi_0$ estimators rely on the strong assumption that there are only
p-values following a uniform distribution on $[0,1]$ in a
neighbourhood of 1. However, Pounds \et \cite{PC} recently observed
the violation of this key assumption. They pointed out that some
p-values associated with induced genes may be artificially sent near
to 1, for example when a one-sided test is performed while the
non-tested alternative is true. To overcome this difficulty, we
propose to estimate the density of p-values by some non-regular
histograms, providing a new estimator of $\pi_0$ that remains
reliable in the Pounds' framework
thanks to a relaxed "identifiability assumption".\\
In the context of density estimation with the quadratic loss and
histograms, asymptotic considerations have been used by Scott
(\cite{Sco79}) for instance. A drawback of this approach relies on
regularity assumptions made on the unknown distribution. Some
AIC-type penalized criteria as in Barron \et \cite{BBM} could be
applied as well. However, such an approach depends on some unknown
constants that have to be calibrated at the price of an intensive
simulation step (see \cite{Leb} in the regression framework). As it
is both regularity-assumption free and computationally cheap, we
address the problem by means of cross-validation, first introduced
in this context by Rudemo (\cite{Rud}). More precisely, the
leave-p-out cross-validation (LPO) is successfully applied following
a strategy exposed in Celisse \et \cite{CR07}. Unlike Schweder and
Spjøtvoll's estimator of $\pi_0$ (\cite{SS82}), ours is fully
adaptive thanks to the LPO-based approach, {\it e.g.} it does not
depend on any
user-specified parameter.\\
The paper is organized as follows. In Section 1, we present a
cross-validation based estimator of $\pi_0$ (denoted by
$\widehat{\pi}_{0}$). Our main assumptions are specified and a
description of the whole $\pi_0$ estimation procedure is given.
Section 2 is devoted to asymptotic results such as consistency of
$\wpi_0$. Then we propose a plug-in multiple testing procedure
(plug-in MTP), based on the same idea as that of Genovese \et
\cite{GW04}. It is compared to the BH-procedure in terms of power
and its asymptotic control of the FDR is derived. Section 3 is
devoted to the assessment of our $\pi_0$ estimation procedure in a
wide range of simulations. A comparison with other existing and
widespread methods is carried out.
The influence of the
$\pi_0$ estimation on the power of the plug-in MTP is inferred as
well. This study results in almost overall improved estimations of
the proposed method.
\section{Estimation of the proportion of true null hypotheses}
\subsection{Mixture model}
Let $P_1,\dots,P_m$ be $m$ \textit{i.i.d.} random variables following
a density $g$ on $[0,1].$ $P_1,\ldots,P_m$ denote the p-values
associated with the $m$ tested hypotheses. Taking into account the
two populations of ($\Hbf_0$ and $\Hbf_1$) hypotheses, we assume
(\cite{Bro,Efr04,GW04}) that $g$ may be written as
\begin{equation*} \forall x\in[0,1],\qquad g(x)=\pi_0
f_0(x)+ (1-\pi_0)f_1(x),
\end{equation*}
where $f_0$ (resp. $f_1$) denotes the density of $\Hbf_0$ (resp.
$\Hbf_1$) p-values, that is p-values corresponding to true null
(resp. false null) hypotheses. $\pi_0$ is the unknown proportion of
true null hypotheses. Moreover, we assume that $f_0$ is continuous,
which ensures that $f_0=1$: $\Hbf_0$ p-values follow the uniforme
distribution $\mathcal{U}([0,1]).$ Subsequently, the above mixture
becomes
\begin{equation}
\forall x\in[0,1],\qquad g(x)=\pi_0 + (1-\pi_0)f_1(x), \label{mixture
model}
\end{equation} where both $\pi_0$ and $f_1$ remain to be estimated.\\
Most of existing $\pi_0$ estimators rely on a sufficient condition
which ensures the identifiability of $\pi_0.$ This assumption may be
expressed as follows
\begin{equation*}
\exists \lambda^*\in]0,1]/\quad \forall i\in\{1,\ldots,m\},\
P_i\in[\lambda^*,1]\Rightarrow P_i \sim
\mathcal{U}([\lambda^*,1]).\qquad(\mathbf{A})
\end{equation*}
(\Abf) is therefore at the origin of Schweder and Spjøtvoll's
estimator (\cite{SS82}), further studied by Storey
(\cite{Sto02,STS}). It depends on a cut-off $\lambda\in[0,1]$ from
which only $\Hbf_0$ p-values are observed. This estimation procedure
is further detailed in Section \ref{simulations}. The same idea
underlies the adaptive Benjamini and Hochberg step-up procedure
described in \cite{BKY}, based on the slope of the cumulative
distribution function of p-values. If we assume $\lambda^*=1$ (that
is $f_1(1)=0$), Grenander \cite{Gre} and Storey \et \cite{ST} choose
$\widehat{g}(1)$ to estimate $\pi_0$, where $\widehat{g}$ denotes
the estimator of $g$. Genovese \et \cite{GW04} use
$(1-G(t))/(1-t),\ t\in(0,1)$ as an upper bound of $\pi_0$, which
becomes (for $t$ large enough) an estimator as soon as
($\Abf$) is true.\\
However, this assumption may be strongly violated as noticed by
Pounds \et\cite{PC}. This point is detailed in Section
\ref{U-shape}. Following this remark, we propose the milder
assumption ($\Abf'$):
\begin{equation*}
\exists \Lambda^*=[\lambda^*,\mu^*]\subset (0,1]/\quad\forall
i\in\{1,\ldots,m\},\ P_i\in\Lambda^*\Rightarrow P_i \sim
\mathcal{U}(\Lambda^*).\qquad(\Abf')
\end{equation*}
While it is a generalization of ($\Abf$), this assumption remains
true in Pounds' framework as we will see in Section \ref{U-shape}\,.
Scheid \et \cite{SS04} proposed a procedure named $Twilight$, which
consists in a penalized criterion and provides, as a by-product, an
estimation of $\pi_0$. Since this procedure does not rely on
assumption (\textbf{A}), it should be taken as a reference
competitor in the simulation study (Section \ref{simulations}) with
respect to our proposed estimators.
\subsection{A leave-$p$-out based density estimator}
\label{LPO density estimation} If $g$ satisfies $(\textbf{A'})$, any
"good estimator" of this density on $\Lambda^*$ would provide an
estimate of $\pi_0.$ Since $g$ is constant on the whole interval
$\Lambda^*,$ we adopt histogram estimators. Note that we do not
really care about the rather poor approximation properties of
histograms outside of $\Lambda^*$ as our goal is essentially the
estimation of $\Lambda^*$ and of the
restriction of $g$ to $\Lambda^*$, denoted by $g_{|\Lambda^*}$ in the
sequel.\\
For a given sample of observations $P_1,\dots,P_m$ and a partition
of $[0,1]$ in $D\in \mathbb{N}^*$ intervals $I=(I_k)_{k=1,\dots,D}$
of respective length $\omega_k=|I_k|,$ the histogram
$\widehat{s}_{\omega}$ is defined by
$$\forall x\in[0,1],\qquad\widehat{s}_\omega(x)=\sum_{k=1}^D
\frac{m_{k}}{m\,\omega_{k}}\ind_{I_{k}}(x),$$ where
$m_k=\sharp\{i\in\lbr1,m\rbr:\,P_i \in I_k\}.$\\
If we denote by $\mathcal{S}$ the collection of histograms we
consider, the "best estimator" among $\mathcal{S}$ is defined in
terms of the quadratic risk:
\begin{eqnarray}
\widetilde{s} & =&
\argmin_{s\in\mathcal{S}}\mathbb{E}_{g}\left[||g-s||_2^2\right],\nonumber\\
 & = &
 \argmin_{s\in\mathcal{S}}\left\{\mathbb{E}_{g}\left[||s||_2^2\right]-2\int_{[0,1]}s(x)g(x)\,dx\right\},
 \label{developped L2-risk}
\end{eqnarray}
where the expectation is taken with respect to the unknown $g$.
According to (\ref{developped L2-risk}), we define $R$ by
\begin{equation}
R(s)=\mathbb{E}_{g}\left[||s||_2^2\right]-2\int_{[0,1]}s(x)g(x)\,dx.
\label{L2-risk}
\end{equation}
In (\ref{L2-risk}) we notice that $R$ still depends on $g$ that is
unknown. To get rid of this, we use a cross-validation estimator of
$R$ that will achieve the best trade-off between bias and variance.
Following (\cite{HTF}), we know that leave-one-out (LOO) estimators
may suffer from some high level variability. For this reason we
prefer the use of leave-$p$-out (LPO), keeping in mind that the
choice of the parameter $p$ will enable the control of the
bias-variance trade-off.\\
At this stage, we refer to Celisse \et\cite{CR07} for an exhaustive
presentation the leave-p-out (LPO) based strategy. Hereafter, we
remind the reader what LPO cross-validation consists in and then,
give the main steps of the reasoning. First of all, it is based on
the same idea as the well-known leave-one-out (see \cite{HTF} for an
introduction) to which it reduces for $p=1.$ For a given $p\in\lbr1,
m-1\rbr$, let split the sample $P_1,\dots,P_m$ into two subsets of
respective size $m-p$ and $p$. The first one, called training set, is
devoted to the computation of the histogram estimator whereas the
second one (the test set) is used to assess the behaviour of the
preceding estimator. These two steps have to be repeated $m\choose p$
times, which is the number of different subsets
of cardinality $p$ among $\{P_1,\dots,P_m\}.$\\
\paragraph{Closed formula of the LPO risk}
This outlined description of the LPO leads to the following closed
formula for the LPO risk estimator of $R(\widehat{s}_\omega)$ (see
\cite{CR07}):
For any partition $I=(I_k)_{k=1,\dots,D}$ of $[0,1]$ in $D$ intervals
of length $\omega_k=|I_k|$ and $p\in\lbr 1,m-1 \rbr,$
\begin{equation}\label{LPO risk}
\widehat{R}_{p}(\omega)=\frac{2m-p}{(m-1)(m-p)}\sum_{k=1}^D\frac{m_{k}}{m\omega_{k}}-\frac{m(m-p+1)}
{(m-1)(m-p)}\sum_{k=1}^D
\frac{1}{\omega_{k}}\left(\frac{m_{k}}{m}\right)^2,
\end{equation} where $m_k=\sharp\{i\in \lbr1,m\rbr:\ P_i\in I_k\},$
$k=1,\dots,D.$
\noindent As it may be evaluated with a computational complexity of
only $O\left(m\log m\right)$, (\ref{LPO risk}) means that we have a
very efficient estimator of the quadratic risk
$R(\widehat{s}_\omega)$.
Now, we propose a strategy for the choice of $p$ that relies on the
minimization of the mean square error criterion (MSE) of our LPO
estimator of the risk. Indeed among
$\{\widehat{R}_{p}(\widehat{s}_\omega):\ p\in \lbr1,m-1\rbr\},$ we
would like to choose the estimator that achieves the best
bias-variance trade-off. This goal is reached by means of the MSE
criterion, defined as the sum of the square bias and the variance of
the LPO risk estimator. Thanks to (\ref{LPO risk}), closed formulas
for both the bias (\ref{bias}) and the variance (\ref{variance}) of
LPO risk estimator may be derived. We recall here these expressions
that come from \cite{CR07}.
\paragraph{Bias and variance of the LPO risk estimator}
Let $\omega$ correspond to a $D-$partition $(I_{k})_{k}$ of $[0,1]$
and for any $k\in \{1,\dots,D\},$ $\alpha_{k}=\Pr[P_{1}\in I_{k}]$
such that $\alpha=(\alpha_{1},\dots,\alpha_{D})\in[0,1]^D.$\\
Then for any $p \in \lbr 1,m-1\rbr$,
\begin{eqnarray}
B_{p}(\omega)&=&B_{p}(\alpha,\omega)=\frac{p}{m(m-p)}\sum_{k=1}^{D}\frac{\alpha_{k}
(1-\alpha_{k})}{\omega_{k}}\,,\label{bias}\\
V_{p}(\omega)&=&V_{p}(\alpha,\omega)=\frac{p^2\varphi_{2}(m,\alpha,\omega)+p\,\varphi_{1}
(m,\alpha,\omega)+\varphi_{0}(m,\alpha,\omega)}{[m(m-1)(m-p)]^2}\,,
\label{variance}
\end{eqnarray}where
\begin{eqnarray*}
\forall (i,j)&\in&\{1,\dots,3\}\times \{1,2\},\quad s_{i,j}=\sum_{k=1}^D\alpha_{k}^i/\omega_{k}^j,\\
\varphi_{2}(m,\alpha,\omega)&=&2m(m-1)\left[(m-2)(s_{2,1}+s_{1,1}-s_{3,2})-m
s_{2,2}-(2m-3)s_{2,1}^2\right]\,,\\
\varphi_{1}(m,\alpha,\omega)&=&-2m(m-1)(3m+1)\left[(m-2)(s_{2,1}-
s_{3,2})-m
s_{2,2}\right]+\\
& &2m(m-1)\left[2(m+1)(2m-3)s_{2,1}^2+
(-3m^2+3m+4)s_{1,1}\right]\,,\\
\varphi_{0}(m,\alpha,\omega)&=&4m(m-1)(m+1)\left[(m-2)(s_{2,1}-s_{3,2})-
m s_{2,2}\right]-\\
& &2m(m-1)\left[(m^2+2m+1)(2m-3)s_{2,1}^2 +(2m^3-4m-2)s_{1,1}\right]+\\
& &m(m-1)^2(s_{1,2}-s_{1,1}^2)\,.
\end{eqnarray*}
\noindent Plug-in estimators may be obtained from the preceding
quantities by just replacing $\alpha_k$ with
$\widehat{\alpha}_k=m_k/m$ in the expressions.
\noindent Following our idea about the choice of $p$, we define for
each (partition) $\omega$ the best theoretical value $p^*$ as the
minimum location of the MSE criterion:
\begin{equation}\label{optimal p parameter}
p^*=\argmin_{p\in\lbr1,m-1\rbr}MSE(p)=\argmin_p\left\{B_{p}(\omega)^2+V_{p}(\omega)\right\}.
\end{equation}
The main point is that this minimization problem has an explicit
solution named $p_{\mathbb{R}}^*$, as stated by Theorem 3.1 in
\cite{CR07}. For the sake of clarity, we recall the MSE expression:
\paragraph{Minimum location expression}
With the same notations as for the bias and the variance, we obtain
for any $x\in\mathbb{R}$,
\begin{equation}
MSE(x)=\frac{x^2[\varphi_{3}(m,\alpha,\omega)+
\varphi_{2}(m,\alpha,\omega)]+x\,\varphi_{1}(m,\alpha,\omega)
+\varphi_{0}(m,\alpha,\omega)}{\left[m(m-1)(m-x)\right]^2}\,,
\nonumber
\end{equation} where\
$\varphi_{3}(m,\alpha,\omega)=(m-1)^2(s_{1,1}-s_{2,1})^2.$\\
\noindent Thus, we define our best choice $\widehat{p}$ for the
parameter $p$ by
\begin{equation} \label{definition p}
\widehat{p}= \left| \begin{array}{ll}
k\left(\widehat{p}_{\mathbb{R}}\right)\,,& \mathrm{if}\
\widehat{p}_{\mathbb{R}}\in [1,m-1]\\
1,& \mathrm{otherwise}
\end{array} \right.,
\end{equation}where $k(x)$ denotes the closest integer near to $x$ and
$\widehat{p}_{\mathbb{R}}$ has the same definition as
$p^*_\mathbb{R}$, but with $\widehat{\alpha}$ instead of $\alpha$ in
the expression.\\
\textit{Remark:} There may be a real interest in choosing adaptively
the parameter $p$, rather than fixing $p=1$. Indeed in the
regression framework for instance, Shao \cite{Shao93} and Yang
\cite{Yang07} underline that the simple and widespread LOO may be
sub-optimal with respect to LPO with a larger $p$. In the linear
regression set-up, Shao even shows that $p/m\to1$ as $m\to+\infty$
is necessary to get consistency in selection.
\subsection{Estimation procedure of $\pi_0$}
\label{estimation procedure}
\subsubsection{Collection of non-regular histograms}
\begin{figure}[h!]
\centering
\makebox{\includegraphics[width=9cm,height=5.5cm]{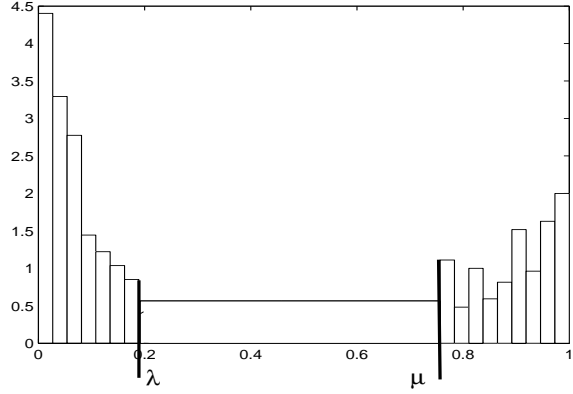}}
\caption{\label{non-regular histogram} Example of non-regular
histogram in $\mathcal{S}$. There are $k=7$ regular columns from $0$
to $\lambda=k/N$, a wide central column from $\lambda$ to
$\mu=\ell/N$, and $N-\ell=7$ regular column of width $1/N$ from
$\mu$ to 1.}
\end{figure}

We now precise the specific collection of histograms we will
consider. For given integers $N_{min}<N_{max},$ we build a regular
grid of $[0,1]$ in $N$ intervals (of length $1/N$) with $N\in\lbr
N_{min}, N_{max}\rbr.$ For a couple of integers $0\leq k<\ell \leq
N$, we define a unique histogram made of first $k$ regular columns of
width $1/N$, then a wide central column of length $(\ell-k)/N$ and
finally $N-\ell$ thin regular columns of width $1/N.$ An example of
such an histogram is given in Figure \ref{non-regular histogram}.
The collection $\mathcal{S}$ of the histograms we consider is defined
by
$$ \mathcal{S}  =  \bigcup_{N\in\lbr
N_{min},N_{max}\rbr}\mathcal{S}_{N},$$ where
\begin{equation*}
\forall N,\qquad \mathcal{S}_N  =  \left\{\widehat{s}_\omega\ :\
 w_{k+1}=(\ell-k)/N,\ w_i=1/N\ for\ i\neq k+1,\ 0\leq
k<\ell\leq N\right\}.
\end{equation*}
Provided $(\mathbf{A'})$ is fulfilled, we expect for each $N$ a
selected histogram with its wide central interval $[\lambda,\mu]$
close to $\Lambda^*$. The comparison of all these histograms (one
per value of $N$) enables to relax the dependence of each selected
histogram on the grid width $1/N.$
\subsubsection{Estimation procedure}\label{algorithm}
Following the idea at the beginning of Section \ref{LPO density
estimation}, $\widehat{\pi}_0$ will consist of the height of the
selected histogram on its central interval $[\lambda,\mu]$. More
precisely, we propose the following estimation procedure for
$\pi_0.$ For each partition (represented here by the vector
$\omega$), we compute
$\widehat{p}\,(\omega)=\argmin_p\widehat{MSE}(p,\omega),$ where
$\widehat{MSE}$ denotes the MSE estimator obtained by plugging
$m_k/m$ in place of $\alpha_k$ in expressions of (\ref{optimal p
parameter}). The best (in terms of the bias-variance trade-off) LPO
estimator of the quadratic risk $R(\widehat{s}_{\omega})$ is
therefore $\widehat{R}_{\widehat{p}(\omega)}(\omega)$. Then we
choose the histogram that reaches the minimum of the latter
criterion over $\mathcal{S}.$ From this histogram, we finally get
both the interval $[\widehat{\lambda},\widehat{\mu}]$, which
estimates $\Lambda^*$, and
$$\widehat{\pi}_0=\widehat{\pi}_0(\widehat{\lambda},\widehat{\mu})\stackrel{def}{=}\frac{\sharp\left\{i:\,P_i \in
\left[\widehat{\lambda},\widehat{\mu}\right]\right\}}{m(\widehat{\mu}-\widehat{\lambda}\,)}\,\cdot$$
These steps are outlined hereafter\\
\textbf{Procedure:}
\begin{enumerate}
\item For each partition denoted by $\omega$,
define\quad $\widehat{p}\,(\omega)=\argmin_p\widehat{MSE}(p,\omega).$
\item Find the best partition\quad
$\widehat{\omega}=\argmin_{\omega}\widehat{R}_{\widehat{p}(\omega)}(\omega).$
\item From $\widehat{\omega}$, get\quad $(\widehat{\lambda},\widehat{\mu}).$
\item Compute the estimator \quad
$\widehat{\pi}_{0}=\frac{\sharp\left\{i:\,P_i \in
\left[\widehat{\lambda},\widehat{\mu}\right]\right\}}{m(\widehat{\mu}-\widehat{\lambda}\,)}\,\cdot$
\end{enumerate}
\section{Asymptotic results}
\subsection{Pointwise convergence of LPO risk estimator}
\begin{lem} \label{MSE convergence} Following the notations in
Section \ref{LPO density estimation}, for any $p\in\{1,m-1\}$ and
$\omega$, we have
\begin{equation*}
MSE(p,\omega)=\mathcal{O}_{m\to+\infty}\left(1/m\right),
\end{equation*}
Moreover if\quad $s_{2,2}+s_{3,2}-s_{1,1}-s_{2,1} \neq 0,$\quad
\begin{equation*}
\widehat{p}(\omega)/m\xrightarrow[m\to+\infty]{a.s.}\ell_{\infty}(\omega),
\end{equation*}
where $\ell_{\infty}(\omega)\in [0,1].$
\end{lem}
\begin{proof}
\begin{enumerate}\item[]
\item We see that
\begin{eqnarray*}
\varphi_3+\varphi_2 & = & 2 m^3 [s_{2,1}+s_{1,1}-s_{3,2}-s_{2,2}-2s_{2,1}^2]+o(m^3),\\
\varphi_1 & = & 2 m^4
[3s_{3,2}+3s_{2,2}-3s_{2,1}+4s_{2,1}^2-3s_{1,1}]+o(m^4),\\
\varphi_0 & = & -4m^5[s_{2,1}^2+s_{1,1}]+o(m^5).
\end{eqnarray*}
Thus for any $p\in\{1,\dots,m-1\}$ and partition of size vector
$\omega\in [0,1]^D$ we have
\begin{equation*}
MSE(p,\omega) =
\mathcal{O}_{m\to+\infty}\left(\frac{1}{m}\right)\,\cdot
\end{equation*}
\item Simple calculations lead to
\begin{equation*}
\frac{p_{\mathbb{R}}^*(\omega)}{m}\xrightarrow[m\to+\infty]{}
\frac{3(s_{2,1}-s_{3,2}-s_{2,2})+7s_{1,1}}{s_{1,1}+s_{2,1}-s_{2,2}-s_{3,2}}=\ell(\alpha,\omega).
\end{equation*}
As for any $k$
$\widehat{\alpha}_k\xrightarrow[m\to+\infty]{a.s.}\alpha_k,$ the
continuous mapping theorem implies the almost surely convergence.
Finally, the result follows by setting
$\widehat{\ell}(\omega)=\ell(\widehat{\alpha},\omega)$ and
$\ell_{\infty}(\omega)=\ind_{\{\ell(\alpha,\omega) \in
[0,1]\}}\ell(\alpha,\omega).$
\end{enumerate}
\end{proof}

\begin{prop}For any given $\omega,$ define $\wp(\omega)$ as in
Section \ref{LPO density estimation} and
$\wL_p(\omega)=\wR_p(\omega)+||g||_2^2.$ If
$\ell_{\infty}(\omega)\neq 1,$ we have
\begin{equation*}
\wL(\omega)\stackrel{def}{=}\wL_{\widehat{p}}(\omega)
\xrightarrow[m\to
+\infty]{P}L(\omega)\stackrel{def}{=}||g-s_{\omega}||_2^2.
\end{equation*}
\end{prop}
\noindent \textit{Remark:} Note that the assumption on
$\ell_{\infty}$ does seem rather natural. It means that the test set
must be (at most) of the same size as the training set
($\widehat{p}/(n-\widehat{p})=\mathcal{O}_{P}(1)$). Moreover,
$\ell_{\infty}(\omega)=1$ if and only if
$s_{2,1}-s_{2,2}-s_{3,2}=-3\,s_{1,1},$ that holds for very specific
densities.
\begin{proof}
The first part of Lemma \ref{MSE convergence} implies that
$\wR_p(\omega)-R(\widehat{s}_\omega)\xrightarrow[m\to +\infty]{P}0.$
Combined with
$R(\widehat{s}_\omega)\xrightarrow[m\to\infty]{}L(\omega)-||g||_2^2,$
it yields that for any fixed $p$,
\begin{equation*}
\wL_p(\omega)\xrightarrow[m\to +\infty]{P}L(\omega).
\end{equation*} Finally, the result follows from both the continuous
mapping theorem and the assumption on $\ell_{\infty}.$
\end{proof}
\subsection{Consistency of $\wpi_0$}
We first emphasize that for a given $N\in\{N_{\min},\dots,N_{\max}\}$
any histogram in $\mathcal{S}_{N}$ is associated with a given
partition of $[0,1]$ that may be uniquely represented by
$(N,\lambda,\mu)$. We give now the first lemma of the consistency
proof.
\begin{lem}\label{minimum location}
For $\lambda^*\neq \mu^*\in[0,1],$ let $s$ be a constant density on
$[\lambda^*,\mu^*]$. Suppose $N_{\min}$ such that for any
$N_{\min}\leq N ,$ it exists a partition $(N,\lambda,\mu)$
satisfying $0<\mu-\lambda\leq \mu^*-\lambda^*.$ For a given $N,$ let
$\omega_N$ represent the partition $(N,\lambda_N,\mu_N)$ with
$\lambda_N=\lceil N\lambda^*\rceil/N$ and $\mu_N=\lfloor
N\mu^*\rfloor/N.$ Define $s_{\omega}$ as the orthogonal projection
of $s$ onto piecewise constant functions built from the partition
associated with $\omega.$ If the dimension of a partition is its
number of pieces, then $\omega_N$ is the partition with the smallest
dimension satisfying
\begin{equation*}
\omega_N \in \argmin_{\omega}||s-s_\omega||_2^2.
\end{equation*}
\end{lem}
\begin{proof}
For symmetry reasons, we deal with partitions, for a given $N,$ made
of regular columns of width $1/N$ from 0 to $\lambda$ and only one
column from $\lambda$ to 1 (\textit{e.g.} we set $\mu=1$). In the
sequel, $I^{(N)}$ denotes the partition associated with $\omega_N$.
\begin{enumerate}
\item Suppose that it exists $\omega_0$ such that $s=s_{\omega_0}.$ Then
$||s-s_{\omega_{N}}||_2^2=0$ and
$\omega_N\in\argmin_{\omega}||s-s_\omega||_2^2$.
\item Otherwise, $s$ does not equal to any $s_\omega$.
\begin{enumerate}
\item If $\lambda^*=k/N,$ then $\lambda_N=\lambda^*.$
Any subdivision $I$ of $I^{(N)}$ satisfies
$||s-s_\omega||_2^2=||s-s_{\omega_N}||_2^2,$ where $\omega$
corresponds to $I$. Now, let $\mathcal{F}_I$ be the set of piecewise
constant functions built from a partition $I$. For any partition
$I=(I_k)_k$ such that $\forall k,\ I^{(N)}_\ell \subset I_k$ for a
given $\ell,$ then $\mathcal{F}_I\subset \mathcal{F}_{I^{(N)}}.$ Thus
$||s-s_\omega||_2^2 =
||s-s_{\omega_N}||_2^2+||s_{\omega_N}-s_\omega||_2^2,$ since
$s_{\omega_N}-s_\omega \in \mathcal{F}_{I^{(N)}}.$ Therefore,
$\omega_N\in\argmin_{\omega}||s-s_\omega||_2^2.$
\item If $\lambda^*\not\in\{1/N,\dots,1\}$.
As before, any subdivision of $I^{(N)}$ will have the same bias,
whereas it is larger for any partition containing $I^{(N)}$. So,
$\omega_N\in\argmin_{\omega}||s-s_\omega||_2^2.$
\end{enumerate}
\end{enumerate}
\end{proof}
\begin{lem}\label{convergence Lhat}
With the same notations as before, we define
$L(\omega)=||s-s_{\omega}||_2^2$. Let $\widehat{L}$ be a random
process indexed by the set of partitions $\Omega$ such that
$\widehat{L}(\omega')\xrightarrow[m\to+\infty]{P}L(\omega'),$ for
any $\omega'\in\Omega$. If
$\widehat{\omega}\in\argmin_{\omega}\widehat{L}(\omega),$ then
\begin{equation*}
\widehat{L}(\widehat{\omega})\xrightarrow[m\to+\infty]{P}\min\{L(\omega):\,
\omega\in\Omega\}.
\end{equation*}
\end{lem}
\begin{proof}
Set $\Gamma\subset \Omega$ such that $\forall \omega \in \Gamma,\
L(\omega)=\min_{\omega'\in \Omega}L(\omega')$ and define
\sloppy$\delta=\min_{\omega\neq \omega'\in\Gamma}|L(\omega)-L(\omega')|/2$.
For $|\Omega|=k$ and $|\Gamma|=\ell,$ we have the ordered quantities
$L(\omega^1)=\dots=L(\omega^k)<L(\omega^{k+1})\leq \dots\leq
L(\omega^\ell)$. Set $\epsilon>0$. For each $\omega^i,$ it exists
$m_i$ (large enough) such that for $m\geq m_i,\
|\wL(\omega^i)-L(\omega^i)|< \epsilon$, with high probability. For
$m_{\max}=\max_i m_i$, we get
$\max_{\omega\in\Omega}|\wL(\omega)-L(\omega)| < \epsilon$ in
probability. Thanks to the latter inequality and by definition of
$\womega$,
\begin{eqnarray*}
L(\womega) < \wL(\womega) + \epsilon \leq \wL(\omega)+ \epsilon <
L(\omega)+ 2\epsilon, \ \mathrm{in\ Probability}
\end{eqnarray*} for any $\omega \in \Omega\setminus \Gamma.$ Hence,
we obtain
\begin{eqnarray*}
L(\womega) < \min_{\omega \in\Omega\setminus \Gamma}
L(\omega)=L(\omega^{k+1}),\ \mathrm{in\ Probability}.
\end{eqnarray*}
Thus, $\womega \in \Gamma$ with high probability and the result
follows.
\end{proof}
\begin{thm}For $0\leq \lambda^*<\mu^*\leq 1,$
let $s:\;[0,1]\mapsto[0,1]$ be a constant function on
$[\lambda^*,\mu^*]$ such that $s$ is not constant on any interval
$I$ with $[\lambda^*,\mu^*]\varsubsetneq I$ (if it exists). Suppose
$N_{\min}$ such that for any $N_{\min}\leq N \leq N_{\max},$ it
exists a partition $(N,\lambda,\mu)$ satisfying $0<\mu-\lambda\leq
\mu^*-\lambda^*.$ Set $\Omega=\cup_{N}\Omega_N$, where $\Omega_N$
denotes the partitions associated with $S_N$. If $\wpi_0$ is the
estimator described in Section \ref{algorithm} selected from
$\Omega$, then
$$\wpi_0\xrightarrow[m\to+\infty]{P} \pi_0.$$
\end{thm}
\begin{proof}
For $\epsilon>0$ and $N_{\min}\leq N \leq N_{\max}$,
\begin{eqnarray*}
\Pr\left[\,|\pi_0-\wpi_0|>\epsilon\right] & = &
\Pr\left[\left|\:s\left(\frac{\lambda^*+\mu^*}{2}\right)-
\widehat{s}_{\widehat{\omega}}\left(\frac{\widehat{\lambda}+\widehat{\mu}}{2}\right)\right|>\epsilon
\right],\\
& \leq & \Pr\left[\:[\widehat{\lambda},\widehat{\mu}]\not\subset
[\lambda^*,\mu^*]\right]+ \Pr\left[\left||s_{\widehat{\omega}}-
\widehat{s}_{\widehat{\omega}}\right||_{2,[\widehat{\lambda},\widehat{\mu}]}^2>
\epsilon^2 (\widehat{\mu}-\widehat{\lambda})\right],\\
& \leq & \Pr\left[\:|L(\widehat{\omega})-L(\omega_N)|>\delta\right]+
\Pr\left[\sup_{\omega}\left||s_\omega-\widehat{s}_{\omega}\right||_2^2>
\epsilon^2 /N_{\max}\right],\\
\end{eqnarray*} for some $\delta>0$ ($||\cdot||_{2,[\widehat{\lambda},\widehat{\mu}]}$ denotes the
quadratic norm restricted to $[\widehat{\lambda},\widehat{\mu}]$).
As the cardinality of the set of partitions is finite ($N_{\max}$
does not depend on $m$),
\begin{eqnarray*}
\Pr\left[\sup_{\omega}\left||s_\omega-\widehat{s}_{\omega}\right||
_2^2>\epsilon^2 /N_{\max}\right] \xrightarrow[m\to+\infty]{} 0.
\end{eqnarray*}
We use the following inequality
$|L(\widehat{\omega})-L(\omega_{N})|-
|L(\widehat{\omega})-\widehat{L}(\widehat{\omega})| \leq
|\widehat{L}(\widehat{\omega})-L(\omega_{N})|$ and the uniform
convergence in probability of $\wL-L$ over $\Omega$
($|\Omega|<+\infty$) to get
\begin{eqnarray*}
\Pr\left[\:|L(\widehat{\omega})-L(\omega_N)|>\delta\right] & \leq &
\Pr\left[\:|\widehat{L}(\widehat{\omega})-L(\omega_N)|>\delta'\right],
\end{eqnarray*} for some $\delta'>0.$
The result comes from both Lemma \ref{minimum location} and Lemma
\ref{convergence Lhat}.
\end{proof}

\subsection{Asymptotic optimality of the plug-in MTP}
\label{MTP procedure} The following is inspired by both \cite{GW04}
and \cite{STS}. In the sequel, we will remind some of their results
to state the link. First of all for any $\theta\in[0,1]$, set
\begin{eqnarray*}
\forall t\in(0,1],\qquad Q_{\theta}(t)  = \frac{\theta\,t}{G(t)}&\
\mathrm{and}\ & \widehat{Q}_{\theta}(t) =
\frac{\theta\,t}{\widehat{G}(t)},
\end{eqnarray*} where $G$ (resp. $\widehat{G}$) denotes the (empirical) cumulative
distribution function of p-values. Let define the threshold
$T_{\alpha}(\theta)=T(\alpha, \theta,
\widehat{G})=\sup\{t\in(0,1):\, \widehat{Q}_{\theta}(t)\leq
\alpha\}.$ Now we are in position to define our plug-in procedure:
\begin{defn}[Plug-in MTP]
Reject all hypotheses with p-values less than or equal to the
threshold $T_\alpha(\wpi_0).$
\end{defn}
\noindent Storey \et \cite{STS} established the equivalence between
the BH-procedure and the procedure consisting in rejecting
hypotheses associated with p-values less than or equal to the
threshold $T_\alpha(1)$, named the step-up $T_\alpha(1)$ procedure.
We may slightly extend Lemma 1 and Lemma 2 in \cite{STS} by using
similar proofs, so that they are omitted here.
\begin{lem}
With the same notations as before, we have
\begin{enumerate}
\item[(i) ] the step-up procedure $T_\alpha(\wpi_0(0,1))=T_\alpha(1)$
is equivalent to the BH-procedure in that they both reject the same
hypotheses,
\item[(ii) ] the step-up procedure
$T_\alpha(\wpi_0(\widehat{\lambda},\widehat{\mu}))$ is equivalent to
the BH-procedure with $m$ replaced by
$\wpi_0(\widehat{\lambda},\widehat{\mu})$.
\end{enumerate}
\end{lem}
\noindent Thus, we observe that the introduction of $\wpi_0$
(supplementary information) in our procedure entails the rejection
of at least as much hypotheses as the BH-procedure ($T_{\alpha}$ in
nonincreasing). Hence our plug-in procedure should be more powerful,
provided it controls the FDR at the required
level $\alpha$.\\
We settle this question now, at least asymptotically, thanks to a
slight generalization of Theorem 5.2 in \cite{GW04} to the case where
$G$ is not necessarily concave (see the "U-shape" framework described
in Section \ref{U-shape} for instance). For $t\in[0,1],$ let define
$FP(t)$ (resp. $R(t)$) as the number of $\Hbf_0$ (resp. the total
number of) p-values lower than or equal to $t$ and set
$\Gamma(t)=FP(t)/(R(t)\vee 1)$. Thus,
$$\forall t\in[0,1],\quad FDR(t)=\E\left[\Gamma(t)\right].$$
\begin{thm}\label{asymptotic control}
For any $\delta>0$ and $\alpha\in[0,\pi_0[,$ define
$\wpi_0^{\delta}=\wpi_0+\delta$. Assume that the density $f$ of
$\Hbf_1$ p-values is differentiable and is nonincreasing on
$[0,\lambda^*],$ vanishes on $[\lambda^*,\mu^*]$ and is nondecreasing
on $[\mu^*,1].$ Then
\begin{enumerate} \item[(i) ]
$Q_{\pi_0}$ is increasing on $I_{\alpha}=Q_{\pi_0}^{-1}([0,\alpha]),$
\item[(ii) ] $\E\left[\Gamma(T_{\alpha}(\wpi_0^{\delta}))\right]\leq
\alpha +o(1). $
\end{enumerate}
\end{thm}
\noindent \textit{Remarks:}\\
Note that the only interesting choice of $\alpha$ actually lies in
$[0,\pi_0)$. If $\alpha\geq \pi_0$, then $FDR(t)\leq \alpha$ is
satisfied in the non-desirable case where all hypotheses are
rejected.\\
A sufficient condition on $G$ for the increase of $Q_{\pi_0}$, is
that $G$ were continuously differentiable and $G'(t)< G(t)/t, \forall
t\in (0,1].$ Thus, $G$ may be nondecreasing (not necessarily concave)
and $Q_{\pi_0}$ may increase yet.\\
\\
\noindent To prove Theorem \ref{asymptotic control}, we first need a
useful lemma, the technical proof of which is deferred to Appendix.
\begin{lem}\label{useful lemma}
With the above notations, for any $\alpha\in(0,1]$,
$T(\alpha,\cdot,\wG):\;[0,1]\mapsto [0,1]$ is continuous a.s.\,.
Moreover for any $\theta\in[0,1],\ G\mapsto T(\alpha,\theta,G)$ is
continuous on $\mathcal{B}^+([0,1])$, the set of positive bounded
functions on $[0,1]$, endowed with the $\lnorm\cdot||_{\infty}$.
\end{lem}
\begin{proof}(Theorem \ref{asymptotic control})
\begin{enumerate}
\item[(i) ] As $f$ is differentiable and nonincreasing, $G$ is concave on
$[0,\mu^*]$ and $Q_{\pi_0}$ increases on this interval. Following
the above remarks, $Q_{\pi_0}$ is still increasing provided
$G'(t)<G(t)/t$ for $t\in[\mu^*,1].$ Thus provided $G'(t)<G(t)/t,\
\forall t\in[\mu^*,1],$ $Q$ increases on $[\mu^*,1]$. Otherwise,
there exists $t_0\in[\mu^*,1]$ such that $G'(t_0)=G(t_0)/t_0$. Then,
the increase of $f$ ensures that $G(x)/x\leq G'(x),\ \forall x\geq
t_0.$ Hence, $Q_{\pi_0}$ is nonincreasing on $[t_0,1].$ Finally
since
$Q(\pi_0)=1$, $Q_{\pi_0}$ is increasing on $I_{\alpha}$.\\
\item[(ii) ]Rewrite first the difference
\begin{eqnarray}
\Gamma\left(T(\alpha, \wpi_0^{\delta},\wG)\right)-\alpha  & = &
\Gamma\left(T(\alpha, \wpi_0^{\delta},\wG)\right)-Q_{\pi_0}\left(T(\alpha, \wpi_0^{\delta},\wG)\right)\nonumber\\
& & + \,Q_{\pi_0}\left(T(\alpha, \wpi_0^{\delta},\wG)\right) - Q_{\pi_0}\left(T(\alpha, \pi_0^{\delta},\wG)\right)
\label{second term}\\
 & & +\, Q_{\pi_0}\left(T(\alpha, \pi_0^{\delta},\wG)\right) -Q_{\pi_0}\left(T(\alpha, \pi_0^{\delta},G)\right)
 \label{third term}\\
& &  + \,Q_{\pi_0}\left(T(\alpha,
\pi_0^{\delta},G)\right)-\alpha\label{last term}.
\end{eqnarray}

\noindent Set $\eta>0$ such that $2\eta<
T(\alpha,\pi_0^{\delta},G)$. Note that
\begin{eqnarray*}
\Gamma\left(T(\alpha,
\wpi_0^{\delta},\wG)\right)-Q_{\pi_0}\left(T(\alpha,
\wpi_0^{\delta},\wG)\right)\leq \frac{1}{\sqrt{m}}\lnorm
\sqrt{m}\left(\Gamma-Q_{\pi_0}\right)||_{\infty,[\eta,1]}+\ind_{\{T(\alpha,
\wpi_0^{\delta},\wG)\leq \eta\}}.
\end{eqnarray*}
Thus thanks to Lemma \ref{useful lemma},
\begin{eqnarray*}
\P\left[T(\alpha, \wpi_0^{\delta},\wG)\leq \eta\right]& \leq &
\P\left[T(\alpha, \pi_0^{\delta},G)\leq \eta+o_{P}(1)\right]
\xrightarrow[m\to+\infty]{}0.\\
\end{eqnarray*}
Besides, both Theorem 4.4 of \cite{GW04} and Prohorov's theorem
(\cite{VdV}) imply that
\begin{eqnarray*}
\E\left[\frac{1}{\sqrt{m}}\lnorm
\sqrt{m}\left(\Gamma-Q_{\pi_0}\right)||_{\infty,[\eta,1]}\right]=o(1).
\end{eqnarray*}
Hence $\E\left[\Gamma\left(T(\alpha,
\wpi_0^{\delta},\wG)\right)-Q_{\pi_0}\left(T(\alpha,
\wpi_0^{\delta},\wG)\right)\right]=o(1).$\\
Thanks to Lemma \ref{useful lemma}, the uniform continuity of
$Q_{\pi_0}$ combined with the convergence in probability of
$\wpi_0^{\delta}$ ensure that the expectation of (\ref{second term}) is of the
order of $o(1).$\\
Since $T(\alpha, \pi_0^{\delta},G)=\sup\{t:\ Q_{\pi_0}(t)\leq
\alpha\pi_0/\pi_0^{\delta}\},$ $\beta=\pi_0/\pi_0^{\delta}<1$ and
$Q_{\pi_0}$ is a one-to-one mapping on $I$, we get
$Q_{\pi_0}\left(T(\alpha, \pi_0^{\delta},G)\right)
=Q_{\pi_0}\left(Q_{\pi_0}^{-1}(\alpha\beta)\right)= \alpha\beta\,$.
Thus,
\begin{equation*}
Q_{\pi_0}\left(T(\alpha, \pi_0^{\delta},\wG)\right)
-Q_{\pi_0}\left(T(\alpha, \pi_0^{\delta},G)\right)
=Q_{\pi_0}\left(T(\alpha\beta, \pi_0,\wG)\right)- \alpha\beta\,,
\end{equation*} Theorem 5.1 (\cite{GW04}) applied with
$\alpha\beta$ instead of $\alpha$ and
$t_0=Q_{\pi_0}^{-1}(\alpha\beta)$ entails that the expectation of
(\ref{third term}) is $o(1)$ as well.\\
Finally, (\ref{last term}) is equal to $(\beta-1)\alpha<0$.
\end{enumerate}
\end{proof}
\section{Simulations and Discussion}\label{simulations}

\subsection{Comparison in the usual framework ($\mu=1$)}
By "usual framework", we mean that the unknown $f_1$ in the mixture
(\ref{mixture model}) is a decreasing density satisfying assumption
(\textbf{A}): it vanishes on an interval $[\lambda^*,1]$ with
$\lambda^*$ possibly equal to 1. In this framework,
$$\widehat{\pi}_0=\frac{\sharp\{i/\,P_i \in[\widehat{\lambda},1]\}}
{m\,(1-\widehat{\lambda}\,)}\,\cdot$$ Except $\widehat{\lambda}$,
this general expression was introduced by Schweder \et \cite{SS82}.
Their estimator
$$\widehat{\pi}_0^{SS}(\lambda)=\frac{\sharp\{i/\,P_i \in[\lambda,1]\}}
{m\,(1-\lambda\,)},$$ is based on (\textbf{A}) and strongly depends
on the parameter $\lambda\in[0,1]$ that is supposed to be given, but
totally unknown in practice. A crucial issue (\cite{LLF}) is
precisely the determination of an 'optimal' $\lambda$.\\

\subsubsection{A potential gain in choosing $\lambda$}
In 2002, Storey \cite{Sto02} studied further this estimator and even
proposed (\cite{ST}) the systematic value $\lambda=0.5$ as a quite
good choice. In the following, we show that even if assumption
(\textbf{A}) is satisfied for $\lambda^*=0.2$ or $0.4$, there is a
real potential gain in choosing $\lambda$ in an adaptive way.\\
In the following simulations, the unknown density $f_1$ in the
mixture (\ref{mixture model}) is a beta density on $[\lambda^*,1]$
with parameter $s$:
$$f_1(t)= s/\lambda^*(1 - t/\lambda^*)^{s-1} \ind_{[0,\lambda^*]}(t),$$
where $(\lambda^*,s)\in\left\{(0.2,4),(0.4,6)\right\}.$ The beta
distribution is all the more sharp in the neighbourhood of 0 as $s$
is large. The proportion $\pi_0$ is equal to 0.9, the sample size
$m=1000$ while $n=500$ repetitions have been made. There does not
seem to be any strong sensitivity to the choice of $N_{max}$ (data
not shown here), as long as $N_{max}$ is obviously not too small.
Until the end of the
paper, $N_{min}=1$ and $N_{max}=100$.\\
Table \ref{choice of lambda} shows the simulation results for the
leave-$p$-out ($LPO$) and the leave-one-out ($LOO$) based estimators
of $\pi_0$, compared to that of Schweder and Spjøtvoll for
$\lambda=0.5$ denoted by $\widehat{\pi}_0^{St}$. We see that in both
cases, $LPO$ is less biased than $LOO$ but slightly more variable,
which leads to a higher value for the MSE. This larger variability
may be due to the supplementary randomness induced by the choice of
$\widehat{\lambda}$. Both $LPO$ and $LOO$ seem a bit conservative
unlike $\widehat{\pi}_0^{St},$ which is however a little less
biased. We say that an estimator of $\pi_0$ is conservative as soon
as it upperbounds $\pi_0$ on average. The main conclusion is that
the MSE of $LPO$ (and $LOO$) is always lower than that of
$\widehat{\pi}_0^{St}$, even if the assumption (\textbf{A}) is
satisfied ($\lambda=0.5>\lambda^*$). An adaptive choice of $\lambda$
may provide a more accurate estimation of $\pi_0$, which is all the
more important as $m$ grows.

\begin{table}
\caption{\label{choice of lambda}Results for the two simulation
conditions $(\lambda^*,s)=(0.2,4)$ and $(\lambda^*,s)=(0.4,6)$. The
LPO and LOO based methods are compared to the Schweder and Spjøtvoll
estimator, $\widehat{\pi}_0^{St}$ computed with $\lambda=0.5$. (All
displayed quantities are multiplied by 100.)} \centering
\begin{tabular}{|c||c|c|c||c|c|c|}
\hline $\pi_0=0.9$& \multicolumn{3}{|c||}{$\lambda^*=0.2,\ s=4$} &
\multicolumn{3}{|c|}{$\lambda^*=0.4,\ s=6$ }\\
\hline \hline
Method & Bias & Std & MSE & Bias & Std & MSE\\
\hline
 $LPO$ & \textbf{0.39} & 2.5 &  \textbf{6.41\ $10^{-2}$} &  \textbf{0.56}& 2.8 & \textbf{8.00\ $10^{-2}$}\\
 $LOO$ &  \textbf{0.46} & 2.3& 5.52\ $10^{-2}$& \textbf{0.61} & 2.7& 7.66\ $10^{-2}$\\
 $\widehat{\pi}_0^{St}$& -0.15& 3.2&  \textbf{9.94\ $10^{-2}$}&0.24 & 3.1& \textbf{ 9.58\ $10^{-2}$}\\
\hline
\end{tabular}
\end{table}

\subsubsection{Comparison when $\lambda^*=1$}
\label{left side simulations} We consider now the general (more
difficult) case when (\textbf{A}) is only satisfied for
$\lambda^*=1$. Thus, $f_1$ is a beta density of parameter $s:$
$f_1(t)=s(1-t)^{s-1},\quad t\in [0,1]$, with $s\in\{5,10,25,50\}.$
The sample size $m=1000$ and $\pi_0\in\{0.5,0.7,0.9,0.95\}.$ Each
condition has been repeated $n=500$ times. We detail below four of
the different methods that have been compared in this framework.

\paragraph{\textit{Smoother} and \textit{Bootstrap}}
\ \\
\noindent In \cite{ST}, the authors proposed a method consisting in
first computing the Schweder and Spjøtvoll estimator on a regular
grid of $[0,1]$ and then adjusting a cubic spline. The final
estimator of $\pi_0$ is the resulting function evaluated
at 1. This procedure is called $Smoother.$\\
The $Bootstrap$ method was introduced in \cite{STS}. Authors define
the optimal value of $\lambda$ as the minimizer of the MSE of their
$\pi_0$ estimator. Since this quantity is unknown, they use an
estimation based on bootstrap. They also need to compute
$\widehat{\pi}_0(\lambda)$ for values of $\lambda$ on a preliminary
grid of $[0,1]$.\\
These methods are available as options of the \textit{qvalue}
function in the {\tt R-}\,package qvalue \cite{ST}.

\paragraph{Adaptive Benjamini-Hochberg procedure}
\ \\
In the sequel, this procedure is denoted by $ABH$ and we refer to
\cite{BKY} for a detailed description. In outline, the method relies
on the idea that the plot of p-values versus their ranks should be
(nearly) linear for large enough p-values (likely $\mathbf{H}_0$
p-values). The inverse of the resulting slope provides a
plausible estimator based on assumption (\textbf{A}).\\
The $ABH$ procedure may be applied through the function {\it
pval.estimate.eta0} in package fdrtool with the option method=
"adaptive" {\tt
http://cran.r-project.org/src/contrib/\\Descriptions/fdrtool.html}.

\paragraph{\textit{Twilight}}
\ \\
In their article, Scheid \et \cite{SS04} proposed a penalized
criterion based on assumption (\textbf{A'}). This is a sum of the
Kolmogorov-Smirnov score and a penalty term. The whole criterion is
expected to provide the widest possible set of $\mathbf{H}_0$
hypotheses. How the penalty term balances against the
Kolmogorov-Smirnov score depends on a constant $C$ that is to be
determined. To do so, the authors propose to use bootstrap combined
with Wilcoxon tests. Besides, this procedure is iterative and
strongly depends on the length of the data, which could be a serious
drawback
with increasing data sets.\\
The function {\it twilight} is available in package twilight
\cite{SS05}.

\begin{figure}[h!]
\centering
\hspace{-1.25cm}
\makebox{\includegraphics[width=15cm,height=10cm]{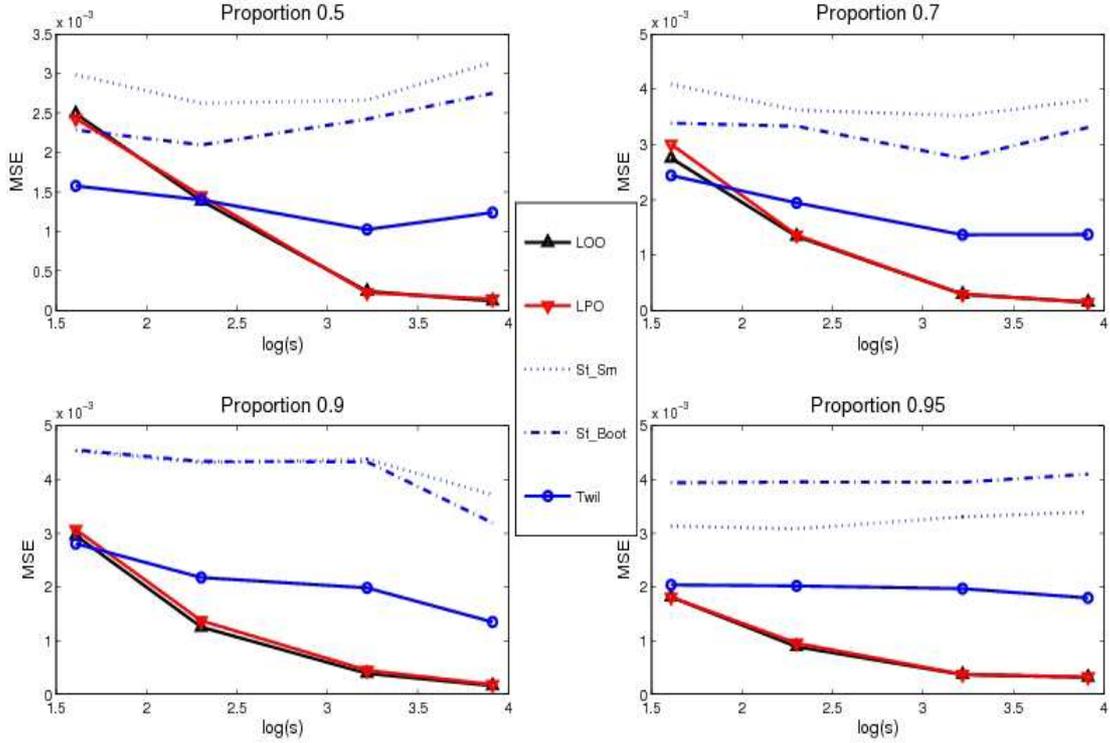}}
\caption{\label{unilateral comparison} Graphs of the MSE of the
$\pi_0$ estimator versus $\log s$, where $s$ is the parameter of the
Beta density. Each graph is devoted to a given proportion, from 0.5
to 0.95\,. $St_{Sm}$ denotes the MSE obtained for $Smoother$,
$St_{Boot}$ that of $Bootstrap$ while $Twil$ states for $Twilight$.}
\end{figure}

\paragraph{Results}
\ \\
As in the preceding simulation study, $LPO$ and $LOO$ refer to the
proposed methods. Figure \ref{unilateral comparison} illustrates the
performances for all the methods but $ABH,$ for which results are
quite poor with respect to other methods (see Table \ref{unilateral
data}). We notice that both $St_{Sm}$ and $St_{Boot}$ have
systematically larger MSE than the three remaining approaches. Our
methods give quite similar results to each other in this framework.
$Twilight$, $LPO$ and $LOO$ furnish nearly the same MSE values in
the most difficult case $s=5$, when $\pi_0>0.5$. Except for
$\pi_0=0.5$ and $s=5,$ $LPO$ and $LOO$ all the more outperform upon
$Twilight$ as the proportion raises. The better performance of
$Twilight$ in this set-up may be due to the classical difference
between cross-validation and penalized criteria. Indeed in the
context of supervised classification for instance, Kearns \et
\cite{KMNR} and Bartlett \et \cite{BBL} show that cross-validation
is used to providing good results, provided the noise level of the
signal is not too high. Otherwise, penalized criteria (like
$Twilight$) outperform upon cross-validation. In the present
context, $s=5$ means that $\Hbf_1$ p-values are spread on a large
part of $[0,1]$ and not only concentrated in a neighbourhood of 0,
while $\pi_0=0.5$ indicates a larger number of $\Hbf_1$ p-values in
the distribution tail of the Beta density. Thus this situation may
be held as the counterpart of the noisy case in supervised
classification. Nevertheless, $LPO$ and $LOO$ always outperform
$Twilight$ when $\pi_0>0.5$. They are even uniformly better than
$Twilight$ for $\pi_0=0.95,$ that is for small proportions of
$\Hbf_{1}$ hypotheses.
\begin{table}
\caption{\label{unilateral data} Numerical results for different
$\pi_0$ estimators with $s=10$ and $\pi_0\in\{0.5, 0.7, 0.9,
0.95\}$. Four other methods are compared to $LPO$ and $LOO$.
$St_{Sm}$ denotes $Smoother$, $St_{Boot}$ states for $Bootstrap$ and
$Twil$ for $Twilight.$(All displayed quantities are multiplied by
100.)} \centering
\begin{tabular}{c}
\begin{tabular}{|c||c|c|c||c|c|c||}
\hline
$\pi_0$ & \multicolumn{3}{|c||}{0.5} &  \multicolumn{3}{|c||}{0.7}  \\
\hline
Method & Bias & Std & MSE  & Bias & Std & MSE \\
\hline \hline
$LPO$ &   1.4  &  3.5   & 14.5\ $10^{-2}$ & 1.4  & 3.4  & 13.6\ $10^{-2}$  \\
\hline
$LOO$  &  1.6   & 3.4  &13.9\ $10^{-2}$  &1.6 &    3.3  & 13.4\ $10^{-2}$  \\
\hline
$St_{Sm}$  & -0.9  & 5.1   & 26.2\ $10^{-2}$  &-0.9  &  6.0 & 36.2\ $10^{-2}$ \\
\hline
$St_{Boot}$  & -2.3  & 4.0 & 20.9\ $10^{-2}$  & -3.3   & 4.7 &  33.3\ $10^{-2}$\\
\hline
$Twil$ &  -1.0  & 3.6 &  14.0\ $10^{-2}$  &-1.5   & 4.2 & 19.4\ $10^{-2}$ \\
\hline
$ABH$ &  37.9   &  8.3  & 15.0    &0.27   & 2.4 & 7.6\ \ \\
\hline
\end{tabular}\\
\\
\begin{tabular}{|c||c|c|c||c|c|c||}
\hline
$\pi_0$ & \multicolumn{3}{|c||}{0.9} & \multicolumn{3}{|c||}{0.95} \\
\hline
Method  & Bias & Std & MSE  & Bias & Std & MSE \\
\hline \hline
$LPO$  & 0.8   & 3.6  & 13.7\ $10^{-2}$  & 0.5  & 3.1  & 9.5\ $10^{-2}$\\
\hline
$LOO$   & 1.0  & 3.4  &  12.5\ $10^{-2}$  & 0.7  & 2.9  & 8.9\ $10^{-2}$\\
\hline
$St_{Sm}$ & -0.5  & 6.6  & 43.1\ $10^{-2}$ & -1.0  & 5.5  & 30.8\ $10^{-2}$\\
\hline
$St_{Boot}$& -3.7  & 5.4  & 43.4\ $10^{-2}$ &    -3.7  & 5.1 &  39.6\ $10^{-2}$\\
\hline
$Twil$  &-1.6  &  4.4  & 21.8\ $10^{-2}$  &-1.6  & 4.2 &  20.2\ $10^{-2}$\\
\hline
$ABH$ & 9.8  & 0.4  &   95.5\ $10^{-2}$ &    4.9   &  0.1  &   24.1\ $10^{-2}$\\
\hline
\end{tabular}
\end{tabular}
\end{table}
\subsection{Comparison in the U-shape case}\label{U-shape}
The 'U-shape case' refers to the phenomenon underlined by Pounds \et
\cite{PC} on a real data set made of Affymetrix 'pooled'
present-absent p-values (one p-value per probe set). 
We explore the behaviour of the preceding methods applied to
p-values with similar distributions. In our simulation design, the
sample is $m=1000$, while $\pi_0\in\{0.25,0.5,0.7,0.8,0.9\}$ and
$n=200$
repetitions of each condition have been made.\\
Typically, the U-shape case appears when one-sided tests are made
whereas the non-tested alternative is true. For example, suppose the
test statistics are distributed as a three-component gaussian mixture
model
\begin{equation}
\pi_0\,\mathcal{N}(0,2.5\:10^{-2})+\frac{1-\pi_0}{2}\left[\mathcal{N}(a,\theta^2)
+\mathcal{N}(b,\nu^2)\right], \label{three component mixture}
\end{equation}
 where $a<0$, $b>0$ and
$\theta,\nu>0,$ corresponding to respectively non-induced,
under-expressed and over-expressed genes. We want to test whether
genes are over-expressed, that is $H_0:$ '{\it the mean equals 0}''
versus $H_1:$ '{\it the mean is positive}'. A test statistic drawn
from $\mathcal{N}(a,\theta^2)$ (under-expressed gene) is more likely
to have a larger p-value than those under $\mathcal{N}(b,\nu^2)$,
which correspond actually to over-expressed genes. This phenomenon
is clearly all the more deep as the gap between $a$ and $b$ is high
and variances $\theta^2$ and $\nu^2$ are small.
Note that a similar shape may be observed when test statistics are
ill-chosen.
\\
In order to mimic Pounds' example, we use (\ref{three component
mixture}) with $-a=b\in\{1,1.5\}$ and $\theta=\nu\in\{0.5,0.75\}.$
As they were quite similar, results in these different conditions
are gathered in Table \ref{bilateral data}.
\begin{table}
\caption{\label{bilateral data} Results of the U-shape case for the
six compared methods for $\pi_0\in\{0.25,0.5,0.7,0.8,0.9\}.$(All
displayed quantities are multiplied by 100.)} \centering
\begin{tabular}{c}
\begin{tabular}{|c||c|c|c||c|c|c||c|c|c||}
\hline $\pi_0$  &  \multicolumn{3}{|c||}{0.25} & \multicolumn{3}{|c||}{0.5} &\multicolumn{3}{|c||}{0.7} \\
\hline \hline
Method &  Bias & Std &  MSE &Bias &Std& MSE &Bias&Std &MSE\\
\hline $LPO$  & 5.5& 6.2& \textbf{0.7}& 5.5 &5.2& \textbf{0.6}& 5.3& 4.4&\textbf{0.5}\\
\hline $LOO$ &6.2&5.7&\textbf{0.7}&6.8&5.7&\textbf{0.8}&6.6&4.8&\textbf{0.7}\\
\hline $St\_{Sm}$ &75.0&0&56.0&50.0&0&25.0&30.0&0&9.0 \\
\hline $St\_{Bo}$ &43.2&3.2&18.7&28.9&2.2&8.4&17.4& 1.6&3.0\\
\hline $Twil$ & 73.2&2.5& 53.6&47.5&3.0&22.6&27.4&2.3&8.0\\
\hline $ABH$ & 45.5&5.4&21.0&31.4&4.2&10.0&19.8&3.1&4.0\\
\hline
\end{tabular}\\
\\
\begin{tabular}{|c||c|c|c||c|c|c||}
\hline $\pi_0$  &\multicolumn{3}{|c||}{0.8} &\multicolumn{3}{|c||}{0.9} \\
\hline \hline Method &   Bias & Std & MSE & Bias & Std & MSE\\
\hline $LPO$  & 5.3 & 4.1 & \textbf{0.4} & 4.2 & 2.7 & \textbf{0.2}\\
\hline $LOO$ & 6.4 & 4.1 & \textbf{0.6}& 4.7 & 2.5 & \textbf{0.3}\\
\hline $St\_{Sm}$ & 20.0 & 0 & 4.0&9.9&0.2&\textbf{1.0} \\
\hline $St\_{Bo}$ & 11.6 &1.3&1.0&5.4&1.6&\textbf{0.3}\\
\hline $Twil$ & 17.5 & 1.8 &3.0&8.0&1.3&\textbf{0.7}\\
\hline $ABH$ & 13.8&2.3&2.0&7.4 &1.3&\textbf{0.6}\\
\hline
\end{tabular}
\end{tabular}
\end{table}
Except $LPO$ and $LOO$ for which this phenomenon is not so strong,
any other method all the more overestimates $\pi_0$ as the
proportion of p-values under the uniform distribution is small. In
our framework, a growth in $\pi_0$ entails an increase in the right
part of the histogram near 1, which is responsible for the
overestimation (violation of assumption (\textbf{A})). On the
contrary when $\pi_0=0.9$, the violation of assumption (\textbf{A}))
is weaker and similar values of MSE are obtained for the competing
approaches. In this set-up, $LPO$, $LOO$ and $St_{Boot}$ provide
systematically the lowest MSE values. In comparison, it is somewhat
surprising that $Twilight$ overestimates $\pi_0$ so much, since it
should have remained reliable under assumption (\textbf{A'}).
Despite the preceding simulation results, we observe a repeated
overestimation, which means that the criterion under-penalizes large
sets of p-values. The involved penalty may have been designed for
the situation before (with only one peak near 0), whereas it may be
no longer relevant in this framework. This may be interpreted as a
consequence of the higher adaptivity of cross-validation based
methods over penalized criteria. Finally it is worth noticing that
both the bias and the MSE of $LPO$ are systematically lower than
those of $LOO$, showing the interest of choosing $p$ in an adaptive
way.
\subsection{Power}
Here, we study the influence of the estimation of $\pi_0$ on the
power of multiple testing procedures obtained as described in
Section \ref{left side simulations} for various $\pi_0$ estimators.
The $Twilight$ method is used for comparison, in association with
the Benjamini-Hochberg procedure (\cite{BH}). Our reference is what
we call the Oracle procedure, which consists in plugging the true
value of $\pi_0$ in the MTP procedure of Section \ref{left side
simulations}. The same simulations as in Section \ref{left side
simulations} are used for this study, which is carried out in two
steps. In the first one, we compare procedures in terms of their
empirical $FDR$, in order to assess the expected control for finite
samples. Thus, we choose the level $\alpha=0.15$ at which we want to
control the $FDR$ and then compute, for each of the $n=500$ samples,
the corresponding $FDP$ in the terminology of \cite{GW04},
\textit{e.g.} the ratio of the number of falsely rejected hypotheses
over the total number of rejections. Finally, we get an estimator of
the actual $FDR$: $\widehat{FDR}$ by averaging the simulation
results. Table \ref{FDR} gives results for the LPO and LOO based
procedures $\widehat{FDR}_{LPO}$, $\widehat{FDR}_{LOO}$ and also for
$Twilight$ ($\widehat{FDR}_{Twil}$), Benjamini-Hochberg
($\widehat{FDR}_{BH}$) and Oracle procedures
($\widehat{FDR}_{Best}$). In the second step, we check the potential
improvement in power enabled by the LPO-based MTP with respect to
the BH-procedure. The assessment of this point is made in terms of
the expectation of the proportion of falsely non-rejected hypotheses
among true alternatives (named $FNR$ here). This criterion is
estimated by the average of the preceding ratio computed from each
sample. Table \ref{power} displays the empirical $FNR$ values,
denoted by $\widehat{FNR}_{LPO}$, $\widehat{FNR}_{LOO}$,
$\widehat{FNR}_{Twil}$, $\widehat{FNR}_{BH}$ and
$\widehat{FNR}_{Best}$ respectively for the LPO, LOO, $Twilight$,
Benjamini-Hochberg and Oracle procedures. In both steps of this
study, $s$ denotes the parameter of the Beta distribution that was
used to simulate the
data.\\
\begin{table}
\caption{\label{FDR} Values of the empirical estimate of the FDR
(\%) for the LPO ($\widehat{FDR}_{LPO}$), LOO
($\widehat{FDR}_{LOO}$), $Twilight$ ($\widehat{FDR}_{Twil}$),
Benjamini-Hochberg ($\widehat{FDR}_{BH}$) and Oracle
($\widehat{FDR}_{Best}$) procedures. $s$ denotes the parameter of
the Beta distribution used to generate the data.}
\centering

\begin{tabular}{|c|c||c|c|c|c|c||}
\hline
s &  $\pi_0$ &  $ \widehat{FDR}_{LPO}$ & $\widehat{FDR}_{LOO}$ &
$ \widehat{FDR}_{Twil}$  & $ \widehat{FDR}_{BH}$ &$ \widehat{FDR}_{Best}$\\
\hline \hline
5&  0.5&14.15& 14.06&14.85&8.35&14.29\\
\cline{2-7}&0.7&14.13&14.03&14.85&10.40&14.50\\
\cline{2-7}&0.9&15.01&15.01&15.73&14.26&14.81\\
\cline{2-7}&0.95&13.23&13.43&13.76&13.13&13.83\\
\hline
10&0.5&14.74&14.69&15.50&6.94&15.02\\
\cline{2-7}&0.7&15.14&15.09&15.61&10.29&15.12\\
\cline{2-7}&0.9&17.91&17.90&18.08&15.85&17.94\\
\cline{2-7}&0.95&14.65&14.65&15.25&14.37&14.95\\
\hline
25&0.5&14.88&14.82&15.51&7.48&15.04\\
\cline{2-7}&0.7&14.69&14.64&15.19&10.47&14.84\\
\cline{2-7}&0.9&15.50&15.57&16.31&13.56&15.92\\
\cline{2-7}&0.95&14.35&14.22&14.51&13.19&14.19\\
\hline
50&0.5&14.76&14.71&15.42&7.40&14.89\\
\cline{2-7}&0.7&14.81&14.77&15.23&10.36&14.87\\
\cline{2-7}&0.9&13.93&13.82&14.79&13.17&13.98\\
\cline{2-7}&0.95&16.12&16.32&16.57&14.65&16.08\\
\hline
\end{tabular}
\end{table}

\begin{table}
\caption{\label{power} Average proportion of falsely non-rejected
hypotheses (\%) for the LPO ($\widehat{FNR}_{LPO}$), LOO
($\widehat{FNR}_{LOO}$), $Twilight$ ($\widehat{FNR}_{Twil}$),
Benjamini-Hochberg ($\widehat{FNR}_{BH}$) and Oracle
($\widehat{FNR}_{Best}$) procedures. $s$ denotes the parameter of
the Beta distribution used to generate the data.}
\centering
\begin{tabular}{|c|c||c|c|c|c|c||}
\hline
s &  $\pi_0$ &  $ \widehat{FNR}_{LPO}$ & $\widehat{FNR}_{LOO}$
& $ \widehat{FNR}_{Twil}$  & $ \widehat{FNR}_{BH}$ &$ \widehat{FNR}_{Best}$\\
\hline \hline
          5&0.5 &93.94&94.22&91.64&99.78&94.16\\
\cline{2-7}&0.7&99.65&99.65&99.59&99.80&99.63\\
\cline{2-7}&0.9 &99.87&99.87&99.86&99.89&99.86\\
\cline{2-7}&0.95&99.91&99.91&99.90&99.92&99.91\\
\hline
         10&0.5 &25.69&25.91&22.01&96.83&23.22\\
\cline{2-7}&0.7 &96.36&96.44&95.08&99.16&96.03\\
\cline{2-7}&0.9 &99.56 &99.56&99.54&99.64&99.56\\
\cline{2-7}&0.95&99.76&99.76&99.76&99.77&99.74\\
\hline
         25&0.5 &0.88&0.90&0.70&17.72&0.79\\
\cline{2-7}&0.7 &22.83&23.04&20.85&61.00&21.93\\
\cline{2-7}&0.9 &97.89&97.89&97.68&98.49&97.86\\
\cline{2-7}&0.95&99.16&99.16&99.06&99.23&99.14\\
\hline
         50&0.5 &0.96&0.92&0.64&1.58&0.72\\
\cline{2-7}&0.7 &2.26&2.30&2.01&10.07&2.19\\
\cline{2-7}&0.9 &82.40&82.47&80.39&88.05&82.08\\
\cline{2-7}&0.95&96.74&96.76&96.60&97.15&96.74\\
\hline
\end{tabular}
\end{table}
\noindent In comparison to the Oracle procedure (with the true
$\pi_0$), Table \ref{FDR} shows that the LPO procedure provides an
actual value of the FDR that is almost always very close to the best
possible one. Moreover in nearly all conditions, LPO outperforms its
LOO counterpart and remains a little bit conservative, {\it e.g.} it
furnishes a FDR that is lower or equal to the desired level
$\alpha$. This observation empirically confirms the result stated in
Theorem \ref{asymptotic control}. Besides as expected, the
estimation of $\pi_0$ entails a tighter control than that of the
BH-procedure where $\widehat{\pi}_0=1$. Unlike the proposed methods,
$Twilight$ fails in controlling the FDR at the desired level since
$\widehat{FDR}_{Twil}$ is very often larger than
$\widehat{FDR}_{Best}$ (the best reachable value), and even larger
than $\alpha$. Subsequently, $Twilight$ should not enter in the
comparison of
methods in terms of power.\\
Table \ref{power} enlightens that proportions of false negatives may
be very high in most of the simulation conditions, as shown by the
Oracle procedure. Nevertheless, $\widehat{FNR}_{LPO}$ remains very
close to the ideal one. As a remark, note that the $Twilight$ FNR
estimates are also close to the Oracle values, but nearly always
lower. As suggested by FDR results, LOO is less powerful that LPO,
whereas both of them outperform by far the BH-procedure. Note that
the proportion of false negatives strongly decreases when $s$ grows,
which means that $\Hbf_1$ p-values are more and more concentrated in
the neighbourhood of 0. As the interval on which assumption
(\textbf{A}) is satisfied is wider, the problem becomes easier.
Besides, we observe a fall in power when $\pi_0$ grows in general.
Indeed for small proportion of true alternatives, the "border"
between the two populations of p-values is more difficult to define
as a large number of $\Hbf_1$ p-values behave like $\Hbf_0$ ones.
Finally note that very often, the LPO procedure shares (nearly) the
same power as the Oracle one.

\subsection{Discussion}
In this article, we propose a new estimator of the unknown
proportion of true null hypotheses $\pi_0$. It relies on first the
estimation of the common density of p-values by use of non-regular
histograms of a special type, and secondly on the leave-$p$-out
cross-validation. The resulting estimator enables more flexibility
than numerous existing ones, since at least it is still convenient
in the "U-shape" case, without any supplementary computational cost.\\
Our estimator may be linked with that of Schweder and Spjøtvoll for
which almost only theoretical results with $\lambda$ fixed have been
obtained by Storey. However unlike the latter, we provide a fully
adaptive procedure that does not depend on any user-specified
parameter. Thus, asymptotic optimality results are here derived with
$\lambda=\widehat{\lambda}$. They assert, for instance, that the
asymptotic exact control of the FDR with our plug-in MTP is reached.\\
Eventually, a wide range of simulations enlighten that the proposed
$\pi_0$ estimator realizes the best bias-variance tradeoff among all
tested estimates. Moreover, the proposed plug-in procedure is
(empirically) shown to provide the expected control on the FDR (for
finite samples), while being a little more powerful than its LOO
counterpart. Moreover, the results in Section \ref{U-shape} confirm
the interest in choosing adaptively the parameter $p$ rather than
the usual $p=1$ value. The LPO procedure is very often almost as
powerful as the best possible one of this type, obtained when
$\pi_0$ is known.

\section{Appendix}
\begin{proof}(Lemma \ref{useful lemma})\\
First, we show that $T(\alpha,\cdot,\wG)$ is right (resp. left)
continuous on $[0,1)$ (resp. $(0,1]$). As it is a similar reasoning,
we only deal with right continuity.\\
Let $(\epsilon_n)_n\in\left(\mathbb{R}_+^*\right)^{\mathbb{N}^*}$
denote a sequence decreasing towards 0. For any $\theta \in(0,1],$
set $\forall n,\ r_n=T(\alpha,\theta+\epsilon_n,\wG)\ a.s.\,.$ Then
$(r_n)_n$ is an almost surely convergent increasing sequence, upper
bounded by $T(\alpha,\theta,\wG)$. To prove that
$T(\alpha,\theta,\wG)$ is its limit, we show that for any
$\delta>0,$ there exists $\epsilon>0$ satisfying
$T(\alpha,\theta+\epsilon,\wG)\geq T(\alpha,\theta,\wG)-\delta.$
Notice that there exists $\eta>0$ s.t.
$T:=T(\alpha,\theta,\wG)=\sup\{t\in[\eta,1]:\,
\widehat{Q}_{\theta}(t)\leq \alpha\}$. Then for $0<\delta<\eta,$ $
T-\delta=\sup\left\{u\in[\eta-\delta,1-\delta]:\
\frac{\theta(u+\delta)}{\wG(u+\delta)}\leq \alpha\right\}. $
Provided $\delta$ is small enough, $\wG(u+\delta)=\wG(u),\ \forall
u.$ Hence, $ T-\delta=\sup\left\{u\in[\eta-\delta,1-\delta]:\
\frac{\theta u}{\wG(u)}+\frac{\theta\delta}{\wG(u)}\leq
\alpha\right\}, $ and $
T(\alpha,\theta+\epsilon,\wG)=\sup\left\{t\in[0,1]:\ \frac{\theta t
}{\wG(t)}+\frac{\epsilon t}{\wG(t)}\leq \alpha\right\}. $ Thus, any
$0<\epsilon< \delta \theta$ provides the result.\\
For the second point, define $G\in\mathcal{B}^+([0,1])$ and for any
sequence $(\epsilon_n)_n\in\left(\mathbb{R}_+^*\right)^{\mathbb{N}}$
decreasing towards 0, let
$(H_{n})_n\in\left(\mathcal{B}^+([0,1])\right)^{\mathbb{N}}$ denote
a sequence of positive bounded functions satisfying $\forall n,\
||G-H_n||_{\infty}\leq \epsilon_n$. Then for large enough $n$, we
have
\begin{eqnarray*}
\frac{\theta t}{G(t)-\epsilon_n}\leq \alpha  &\Leftrightarrow&
\frac{\theta t}{G(t)}\leq \alpha
\left(1-\frac{\epsilon_n}{G(t)}\right),
\end{eqnarray*}
and $\alpha (1-\epsilon_n/||G||_{\infty})\leq \alpha$. Thus,
$r_n=\sup\{t:\ \theta t/(G(t)+\epsilon_n)\leq \alpha\}$ denotes an
increasing sequence that is bounded by $T(\alpha,\theta,G)$.
Moreover as $(\epsilon_n)_n$ decreases towards 0, $r_n$ is as close
as we want to $T(\alpha,\theta,G)$. The same reasoning may be
followed with $r'_n=\sup\{t:\ \theta t/(G(t)-\epsilon_n)\leq
\alpha\},$ which concludes the proof.
\end{proof}

 \nocite{}

\bibliographystyle{plain}

\bibliography{bibliographie}

\end{document}